\documentclass{article}

\usepackage{amsfonts}
\usepackage{graphicx}
\usepackage[toc,page]{appendix}

\usepackage{mathrsfs}

\newtheorem{teo}{Theorem}[section]
\newtheorem{lemma}[teo]{Lemma} 
\newtheorem{cor}[teo]{Corollary} 

\newtheorem{defi}[teo]{Definition}

\begin{document}

\title{Foliation by free boundary constant mean curvature leaves}
\author{Fabio Montenegro}
\maketitle

\begin{abstract}
Let $M$ be a Riemannian manifold of dimension $n+1$  with smooth boundary and $p\in \partial M$. 
We prove that there exists a smooth foliation around $p$ whose leaves are submanifolds of dimension 
$n$, constant mean curvature and its arrive perpendicular to the boundary of M, provided that $p$ is 
a nondegenerate critical point of the mean curvature function of $\partial M$.
\end{abstract}

\section{Introduction}
The strategy of the proof of this result was inspired by \cite{Ye}. In this work, Rugang Ye 
considered the foliation by geodesic spheres around $p\in M$ of small radius and showed that this foliation can be perturbed into a foliation whose leaves are spheres of constant mean curvature, provided that $p$ is a nondegenerate critical point of the scalar curvature function of $M$. So we are going to consider a family of foliations whose leaves are submanifolds 
of $M$ with boundary contained in $\partial M$ and it's arriving perpendicular to the boundary of $M$.
The idea is then to perturb each leaves to obtain, via implicit function theorem, a foliation whose leaves are hemispheres of constant mean curvature
and its arrive perpendicular to the boundary of M, provided that $p$ is 
a nondegenerate critical point of the mean curvature function of $\partial M$.

We refer to \cite{Fernando} for basic terminology in local Riemannian geometry.
Let $(M,g)$ be an (n+1)-dimensional Riemannian manifold with smooth boundary $\partial M$, 
$n\geq 2$. We will denote by $\nabla$ and $\overline{\nabla}$ the covariant derivatives and by 
$R$  and $\bar{R}$ the full Riemannian curvature tensor of $M$ and $\partial M$, respectivily.
The trace of second fundamental form of the boundary will be denoted by $h$. 
We will make use of the index notation for tensors, commas denoting covariant differentiation and we will adopt the summation convention. 
\begin{defi} 
Let $T$ denote the inward unit normal vector field along $\partial M$.
Let $\Sigma$ be a submanifold with boundary $\partial \Sigma$ contained in $\partial M$. 
The unit conormal of $\partial \Sigma$ that points outside $\Sigma$ will be denoted by $\nu$. $\Sigma$ is 
called {\it free boundary} when $\nu = -T$ on $\partial \Sigma$.
\end{defi} 
{\bf Acknowledgements}. Part ot this work was done while the author was visiting the university of Princeton, New Jersey. I would like to thank Fernando C. Marques for his kind invitation and helpful discussions on the subject.

\section{Fermi Coordinate System}
Consider a point $p\in \partial M$ and an orthonormal basis $\{ e_1,\dots ,e_{n}\}$ for 
$T_p\,\partial M$. 
Let $\mathbb{B}_r=\{x\in \mathbb{R}^n\, ;\;|x|<r\}$ be the open ball in $\mathbb{R}^n$.
There are $r_p>0$ and $t_p> 0$ for which we can define the Fermi coordinate system centered at $p$,
$\varphi^0:\mathbb{B}_{r_p}\times [0,t_p)\to M$, given by
\begin{equation}
\label{varphi0}
\varphi^0(x,t)=\mathrm{exp}^M_{\phi^0(x)}(t\,T(x))
\end{equation}
where $\phi^0(x)=\mathrm{exp}_p^{\partial M}(x^ie_i)$ is the normal coordinate system in
$\partial M$ centered at $p$ and $T(x)$ is the inward unit vector normal to the boundary 
at $\phi^0(x)$.

For each $\tau =(\tau_1,\dots ,\tau_{n})\in \mathbb{B}_{r_p/2}$
we will consider the Fermi coordinate system centered at 
$$
c(\tau)=\mathrm{exp}_p^{\partial M}(\tau^ie_i),
$$
which we denote by
$\varphi^{\tau}:\mathbb{B}_{r_p/2}\times [0,t_p)\to M$, and it defined by
\begin{equation}
\label{varphi^tau}
\varphi^{\tau}(x,t)=\mathrm{exp}_{\phi^{\tau}(x)}^M(t\,T(x))
\end{equation}
where 
$$
\phi^{\tau}(x)=\mathrm{exp}_{c(\tau)}\left(x^ie_i^{\tau}\right),
$$ 
$e_i^{\tau}$ are the parallel transport of $e_i$ to $c(\tau)$ along the geodesic 
$c(s\tau)|_{0\leq s\leq 1}$ in $\partial M$,
and $T(x)$ is the inward unit vector normal to the boundary at $\phi^{\tau}(x)\in \partial M$. 

We will denote the metric tensor of $M$ by 
$ds^2$, the coefficients of $ds^2$ in the coordinates system $\varphi^\tau$ by $\mathrm{g}_{ij}^{\tau}(x,t)$,
 and $g^{-1}_{\tau}=(\mathrm{g}_{ij}^{\tau})^{-1} =(\mathrm{g}^{ij}_{\tau})$. 
The expansion of $ g_{ij}^{\tau}$ (up to fourth order) in Fermi coordinates can be found in \cite[ p.1604]{Fernando}.
\begin{lemma}
\label{lemma1}
In Fermi coordinates $(x_1,\dots ,x_{n}, t)$ centered at $c(\tau)\in \partial M$ we have
$g_{\tau}^{tt}(x,t)= 1$, $g_{\tau}^{ti}(x,t)= 0$, and
\begin{equation}
\label{g-1}
\begin{array}{lr}
g^{ij}_{\tau}(x,t)=\delta_{ij}+2h_{ij}t & \hspace{8cm}
\end{array}
\end{equation}
$$
\begin{array}{l}
\hspace{.4cm} +\frac{1}{3}\bar{R}_{ikjl}\,x_kx_l+2h_{ij,\underline{k}}\,t\,x_k+(R_{titj}+3h_{ik}h_{kj})\,t^2  \\ \\
\hspace{.4cm} +\frac{1}{6}\bar{R}_{ikjl,\underline{m}}\,x_kx_lx_m 
+\left(\frac{2}{3}\,\mathrm{Sym}_{ij}(\bar{R}_{ikml}h_{mj})+h_{ij,\underline{kl}}\right)tx_kx_l  \\ \\
\hspace{.4cm} +\left(R_{titj,\underline{k}}+6\,\mathrm{Sym}_{ij}(h_{il,\underline{k}}h_{lj})\right)t^2 x_k \\ \\
\hspace{.4cm} +\left(\frac{1}{3}R_{titj,t}+\frac{8}{3}\,\mathrm{Sym}_{ij}(R_{titk}h_{kj})+4h_{ik}h_{kl}h_{lj} \right)t^3  \\ \\
\hspace{.4cm} +\left(\frac{1}{20}\bar{R}_{ikjl,\underline{mp}}+\frac{1}{15}\bar{R}_{ikql}\bar{R}_{jmqp}\right)x_kx_lx_mx_p \\ \\
\hspace{.4cm} +\left(\frac{1}{3}\,\mathrm{Sym}_{ij}(\bar{R}_{ilpm,\underline{k}}h_{pj})
+\frac{2}{3}\,\mathrm{Sym}_{ij}(\bar{R}_{ikpl}h_{pj,\underline{m}})+\frac{1}{3}h_{ij,\underline{klm}}\right)tx_kx_lx_m   \\  \\
\end{array}
$$
$$
\begin{array}{l}
\hspace{.4cm} +\left(\frac{1}{2}R_{titj,\underline{kl}}+\frac{1}{3}\,\mathrm{Sym}_{ij}(\bar{R}_{ikml}R_{tmtj})
+\frac{7}{3}\,\mathrm{Sym}_{ij}(\bar{R}_{ikml}h_{pj})h_{mp} \right.  \\  \\
\hspace{.4cm} -\frac{4}{3}\,\mathrm{Sym}_{pj}(\bar{R}_{pkml}h_{mj})h_{ip}
-\frac{4}{3}\,\mathrm{Sym}_{ip}(\bar{R}_{ikml}h_{mp})h_{pj} +\frac{4}{3}\bar{R}_{mkpl}h_{im}h_{pj}  \\  \\
\hspace{2.3cm}  \left. +4\,\mathrm{Sym}_{ij}(h_{im,\underline{kl}}h_{mj})-\frac{1}{2}(h_{im}h_{mj})_{,\underline{kl}} 
+4h_{im,\underline{k}}h_{mj,\underline{l}}\right)t^2x_kx_l   \\  \\
\hspace{.4cm} +\left(\frac{1}{3}R_{titj,t\underline{k}}+\frac{8}{3}\,\mathrm{Sym}_{ij}(R_{titl,\underline{k}}h_{lj})
+\frac{8}{3}\,\mathrm{Sym}_{ij}(R_{titl}h_{lj,\underline{k}})\right.  \\  \\  
\hspace{4.5cm}  \left.+8h_{jl,\underline{k}}h_{lm}h_{im}+4h_{ml,\underline{k}}h_{lj}h_{im}\right)t^3x_k  \\  \\
\hspace{.4cm}  +\left(\frac{1}{12}R_{titj,tt}-\frac{1}{3} R_{titk,t}R_{tktj,t}+R_{titk}R_{tktj}+6\,\mathrm{Sym}_{ij}(R_{titk}h_{lj})h_{kl}\right.  \\  \\
\hspace{.4cm}  +\frac{5}{6}\,\mathrm{Sym}_{ij}(R_{titk,t}h_{kj})-\frac{8}{3}\,\mathrm{Sym}_{lj}(R_{tktj}h_{lk})h_{il}
-\frac{8}{3}\,\mathrm{Sym}_{il}(R_{tktl}h_{ik})h_{lj}  \\  \\
\hspace{6cm}  \left.+\frac{13}{3}R_{tktl}h_{ik}h_{lj}+5h_{ik}h_{kl}h_{lm}h_{mj}\right)t^4  \\  \\
\hspace{.4cm} +a^1_{ijlkm}(\tau,x,t)x_ix_jx_kx_lx_m+a^2_{ijlk}(\tau,x,t)tx_ix_jx_kx_l  \\  \\ 
\hspace{.4cm} +a^3_{ijl}(\tau,x,t)t^2x_ix_jx_k+a^4_{ij}(\tau,x,t)t^3x_ix_j+a^5_{i}(\tau,x,t)t^4x_i+a^6(\tau,t)t^5
\end{array}
$$
where every coefficient is computed at $c(\tau)$
and $a^1_{ijlkm}$, $a^2_{ijlk}$, $a^3_{ijl}$, $a^4_{ij}$, $a^5_{i}$, and $a^6$ are smooth
functions. Here underlined indices mean covariant differentiation as tensor on the boundary.
\end{lemma}
{\it Proof}. We write the expansion
$$
g^{\tau}=(g_{ij}^{\tau})=I+C_1+C_2+C_3+C_4+\mathcal{O}(r^5)
$$
where $I$ is the identity matrix and $C_i$ is homogeneous of degree $i$.

For small $r$, we have
$$
g^{-1}_{\tau}=\sum_{j=0}^{\infty}(-1)^j(g^{\tau}-I)^j
$$
and so
$$
g^{-1}_{\tau}=I-C_1+(C_1^2-C_2)+(2C_1C_2-C_3-C_1^3)
$$
$$
+(C_1^4-C_4+C_2^2+2C_1C_3-3C_1^2C_2)+\mathcal{O}(r^5).
$$
The expansion formula (\ref{g-1}) now follows immediately by applying the expansion of $ g_{ij}^{\tau}$
in \cite[ p.1604]{Fernando}. \hfill \fbox

\section{Pertubation by free boundary submanifods}

We will work with the following set of functions
\begin{equation}
\label{C2alpha}
C_T^{2,\alpha}(S_+^{n})=\left\lbrace\varphi\in C^{2,\alpha}(S_+^{n});\; \frac{\partial \varphi}{\partial e_t}=0 \;\; \mathrm{in}\;\;\partial S_+^n \right\rbrace
\end{equation}
where $S^n_+=\{(x,t)\in\mathbb{R}^{n+1};\; t^2+|x|^2=1, \; t\geq 0\}$.

For $\varphi \in C^{2,\alpha}_T(S^n_+)$ we define 
$$
S^{+}_{\varphi}=\{(1+\varphi(x,t))(x,t);\;(x,t)\in S_+^n\}
$$
\begin{figure}[!htp]
\centering 
\includegraphics[width=6cm]{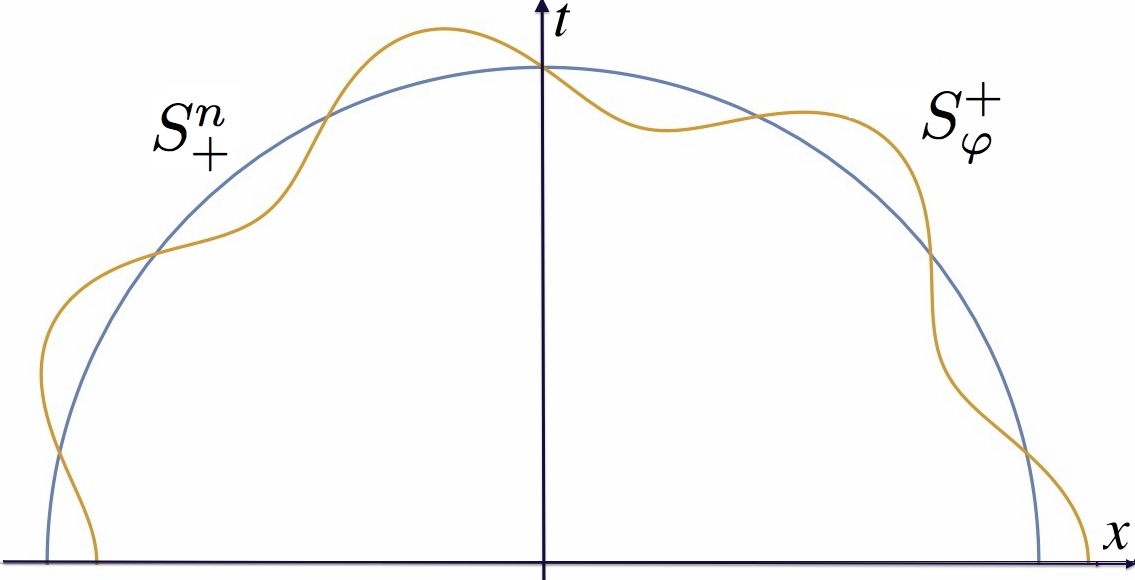} 
\end{figure}
and 
\begin{equation}
\label{Srtauvaphi}
S_{r,\tau ,\varphi}=\varphi^{\tau}(\alpha_r(S^{+}_{\varphi})) 
\end{equation}
where $\alpha_r$ is the dilation $(x,t)\mapsto (rx,rt)$ for $0<r<r_0$ and $r_0$ sufficiently small such that
$$
\alpha_{r_0}(\mathbb{B}_{2}\times [0,2))\subset \mathbb{B}_{r_p/2}\times [0,t_p).
$$

There are numbers $\delta_0>0$ and $r_0>0$ such that $S_{r,\tau ,\varphi}$ is an embedded $C^2$ hypersurface in $M^{n+1}$ for any $\| \varphi \|_{C^1}\leq \delta_0$ and $0<r<r_0$. In addition 
$S_{r,\tau ,\varphi}$ is a free boundary submanifold of $M$, this is, $\partial S_{r,\tau ,\varphi}\subset \partial M$ and
 its arrive perpendicular to the boundary of M, because $\partial \varphi / \partial e_t=0$.
We denote the inward mean curvature function of $S_{r,\tau ,\varphi}$ by $h(r,\tau ,\varphi)$.

For $(x,t)\in S^n_+$ we denote $H(r,\tau ,\varphi)(x,t)$ the inward mean curvature of the surface 
$S_{\varphi}^+$ at $(1+\varphi(x,t))(x,t)$ with respect to the metric 
$ds_{\tau ,r}^2$ on $\mathbb{B}_2\times [0,2)$,
given by $ds_{\tau ,r}^2=r^{-2}\alpha_r^*(\varphi_r^*(ds^2))$, here $ds^2$ is the metric tensor in $M$. 

For each $(x,t)\in S^n_+$ we have
\begin{equation}
\label{H=rh}
H(r,\tau ,\varphi)(x,t)=r\,h(r,\tau ,\varphi)\left(\varphi^{\tau}(r(1+\varphi (x,t))(x,t))\right).
\end{equation}

For $(x,t)\in \mathbb{B}_2\times [0,2)$ and $v,w\in T_{(x,t)}\mathbb{B}_+^2$ we have
$$
ds_{\tau ,r}^2(x,t)(v,w)=\frac{1}{r^2}\varphi_r^*(ds^2)(rx,rt)(d\alpha_r(x,t)v,d\alpha_r(x,t)w)
$$
$$
=\varphi_r^*(ds^2)(rx,tx)(v,w).\hspace{0.45cm}
$$
But, by the Lemma 2.2 in \cite[ p.1604]{Fernando},
$$
ds_{\tau ,r}^2(x,t)(v,w)=\langle v,w\rangle_{\mathbb{R}^{n+1}}+O(r)
$$
with $O(r)\to 0$ when $r=|(x,t)|\to 0$. One readily checks $ds_{\tau ,r}^2$ extends smoothly to the 
euclidean metric when $r$ goes to zero.
Hence $H(r,\tau ,\varphi )$ also extends to $r=0$. Then by a straightforward computation the inward mean curvature function of $S_{\varphi}^+$ at $(\bar{x},\bar{t})=(1+s\varphi (x,t))(x,t)$ with respect to the metric 
$ds_{\tau ,r}^2$ on $\mathbb{B}_2^+$, can be written as

\begin{equation}
\label{Hrtausphi}
H(r,\tau , s\varphi )(\bar{x},\bar{t})=\frac{1}{\Psi_s}\left(\Delta \rho -s\,\Delta \overline{\varphi}
-\frac{s^2}{2}\Delta \overline{\varphi}^2 \right)
\end{equation}
$$
-\frac{1+s\overline{\varphi}}{\Psi_s^2}\left[\frac{\partial \Psi_s}{\partial t}
\left(t-s\,\frac{\partial \overline{\varphi}}{\partial t} \right)
+\sum_{i,j}g_{\tau}^{ij}\,\frac{\partial \Psi_s}{\partial x_i}
\left(x_j-s\,\frac{\partial \overline{\varphi}}{\partial x_j} \right)\right]
$$
where $\rho (x,t) =(t^2+|x|^2)/2$, $g_{\tau}^{ij}=g_{\tau}^{ij}(\bar{x},\bar{t})$,
\begin{equation}
\label{bar(varphi)}
\overline{\varphi}(x,t)=\varphi \left(\frac{x}{\sqrt{t^2+|x|^2}}\,,\frac{t}{\sqrt{t^2+|x|^2}} \right),
\end{equation}
\begin{equation}
\label{Psi-s}
\Psi_s=\Psi_s(r,x,t)=
\end{equation}
$$
\sqrt{\!\left(\!t-s(1\!+\!s\overline{\varphi})\frac{\partial \overline{\varphi}}{\partial t} \right)^2
\!\!\!+\sum_{i,j}g_{\tau}^{ij}(rx,rt)\!\left(\!x_i-s(1\!+\!s\bar{\varphi})\frac{\partial \overline{\varphi}}{\partial x_i} \right)\!\!
\left(\!x_j-s(1\!+\!s\overline{\varphi})\frac{\partial \overline{\varphi}}{\partial x_j} \right)}
$$
and $\Delta$ is the standard Laplace operator on $B_2^+$ relative to the metric $ds_{\tau ,r}^2$.
\begin{lemma}
\label{lemma2}
We have
\begin{equation}
\label{H(r,t,0)}
\begin{array}{lr}
H(r,\tau ,0)(x,t)=n+\left[h_{ii}^{\tau}\,t-(n+3)h_{ij}^{\tau}\,tx_ix_j\right]r 
+\left[\frac{3n+2}{2}h_{ij}^{\tau}h_{kl}^{\tau}\,t^2x_ix_jx_kx_l
\right. & \hspace{1cm}
\end{array}
\end{equation}
$$
\begin{array}{l}
\hspace{.4cm} -(n+4)h_{ij,\underline{k}}^{\tau}\,tx_ix_jx_k
+\left(-\frac{n+4}{2}R_{titj}-\frac{3n+20}{2}\,h_{ik}^{\tau}h_{kj}^{\tau}-h_{ij}^{\tau}h_{kk}^{\tau}\right)t^2x_ix_j \\ \\
\hspace{.4cm} \left. +\frac{1}{3}\bar{R}_{kikj}\,x_ix_j+2h_{ji,\underline{j}}^{\tau}\,tx_i+2(h_{ij}^{\tau})^2t^2 
\right]r^2+\left[\displaystyle\int_0^1\frac{(1-\eta )^2}{2}H_{rrr}(\eta \,r,\tau ,0)\,d\eta\right]r^3
\end{array}
$$
where every coefficient is computed at $c(\tau)$.
\end{lemma}

\begin{cor}
The following holds true
$$
H(0,\tau,0)=\lim_{r\to 0}H(r,\tau ,0)=n.
$$
\end{cor}

Now we consider $H(r,\tau ,\cdot)$ as a mapping from 
$C_N^{2,\alpha}(S_+^{n})$ into $C^{0,\alpha}(S_+^{n})$ 
and let $H_\varphi$ denote the differential of $H$ with respect to $\varphi$. 
In order to calculate $H_\varphi$ we consider the variation of $S_+^n$
by smooth maps $f:S_+^n\times (-\epsilon, \epsilon)\to \mathbb{B}_+^2$ given by 
$f(x,t,s)=(1+s\,\varphi(x,t))(x,t)$. For each $s\in (-\epsilon, \epsilon)$ we
denote $f^s(x,t)= f(x,t,s)$. Note that $f^s(S_+^n)=S_{s\varphi}^+$ is an embedded $C^2$ in 
$\mathbb{B}_+^2$ with $\partial S_{s\varphi}^+\subset \partial \mathbb{B}_+^2$. 
We will denote by $N_s(r,\tau ,\varphi)$ a unit vector field normal to $S_{s\varphi}^+$
and $H(r,\tau ,s\varphi )$ the mean curvature of $S_{s\varphi}^+$. We decompose the variational 
vector field 
$$
\partial_s=\varphi(x,t)(x,t)=\partial_s^T+v^sN_s
$$
where $v^s$ is the function on $S_+^n$ defined by $v^s=ds_{r,\tau}^2(\partial_s,N_s)$.

\begin{figure}[!htp]
\centering 
\includegraphics[width=6.5cm]{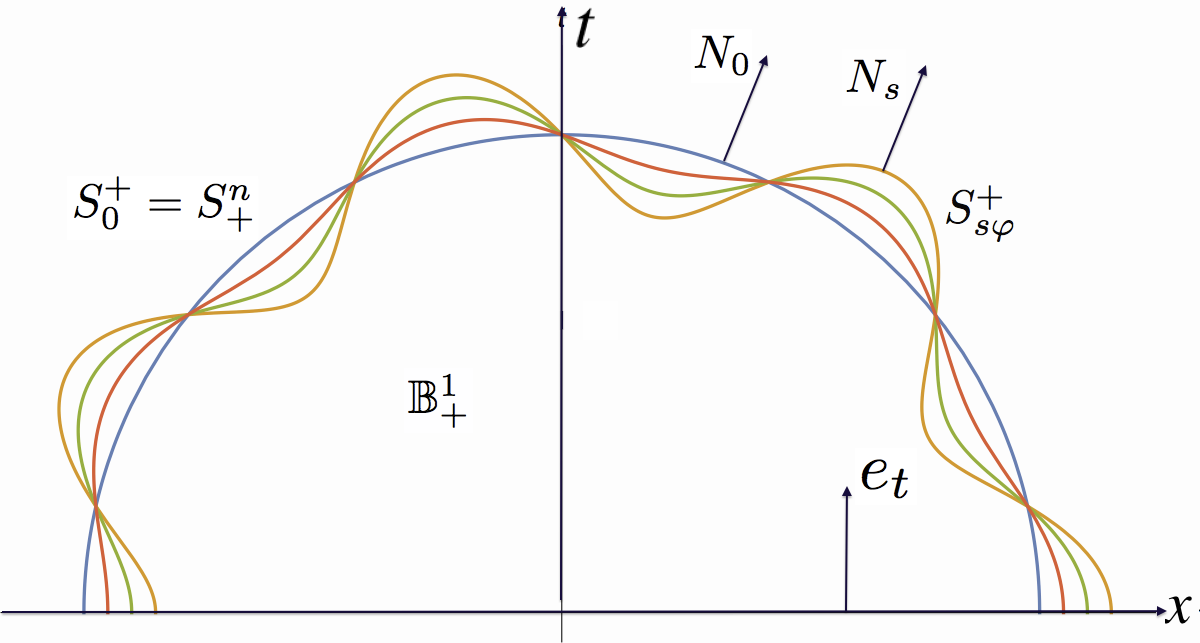} 
\end{figure}

By the Proposition 16 in \cite[ p.14]{Lucas} we have
$$
H_\varphi(r,\tau ,0 )\varphi=(\partial_sH(r,\tau ,s\varphi ))\big\vert_{s=0}=dH(r,\tau ,0)(\partial_0^T)-L_{r,\tau}v_0 
$$
and
$$
\partial_sds_{\tau ,r}^2(N_s,e_t)\big\vert_{s=0}=-\frac{\partial v^0}{\partial e_t}+ds_{\tau ,r}^2(N_0,\nabla_{N_0}e_t)v^0 
$$
where $L_{r,\tau}=\Delta_{(S_{+}^n,ds_{r,\tau}^2)}+Ric_{r,\tau}(N_0,N_0)+\|B_{r,\tau}\|^2$ is the Jacobi operator.

In particular 
\begin{equation}
\label{H_phi}
H_{\varphi}(0,\tau ,0)\varphi =L\varphi:=-(\Delta_{S_+^n}+n)\varphi ,
\end{equation}
where $\Delta_{S_+^n}$ is the standard Laplace operator in $S_+^n$. 
\begin{lemma}
\label{lemma3}
We have
\begin{equation}
\label{H_rphi}
\begin{array}{l}
H_{\varphi r}(0,\tau ,0)\varphi =
2\displaystyle\frac{\partial \overline{\varphi}}{\partial x_i}\,h_{ij}^{\tau}\,x^j\! \left[ (x^{n+1})^3\!+\!(x^{n+1})^2\!+n\,x^{n+1}\right]\!
-2\frac{\partial^2 \overline{\varphi}}{\partial x_i\partial x_j}\,h_{ij}^{\tau}\,x^{n+1}
\end{array}
\end{equation}
$$
+\frac{\partial \overline{\varphi}}{\partial t}\left(h_{ii}^{\tau}+{\partial t}\,h_{ij}^{\tau}\,x^ix^j \right)
+\left( \Delta_{\mathbb{S}^n_+} \varphi +3\varphi \right) h_{ij}^{\tau}\,x^ix^jx^{n+1}
-\varphi\,h_{ij}^{\tau}\,x^ix^j(x^{n+1})^2
$$
and
\begin{equation}
\label{H_phiphi}
H_{\varphi \varphi}(0,\tau ,0)\varphi \varphi =2n\,\varphi^2
-(n-2)\left(\frac{\partial \overline{\varphi}}{\partial t}\right)^2-(n-2)\sum_i\left(\frac{\partial \overline{\varphi}}{\partial x_i}\right)^2
\end{equation}
where $\overline{\varphi}$ was defined in (\ref{bar(varphi)}).
\end{lemma}

The Jacob operator 
$$
L:C_T^{2,\alpha}(S_+^n)\to C^{0,\alpha}(S_+^n)
$$
$L=\Delta_{S_+^n}+n$ has an $n$-dimensional kernel $K$ consisting of first order spherical harmonic functions 
$x^i=x^i|_{S_+^n}$, $i=1,\dots ,n$, which satisfy
$$
\frac{\partial }{\partial e_t}x^i\bigg\vert_{S_+^n}=0\;\; \mathrm{in}\;\;\partial S_+^n.
$$
In addition we have the $L_2$-decompositions of spaces $C_N^{2,\alpha}(S_+^n)=K\oplus K^{\perp}$ and 
$C^{0,\alpha}(S_+^n)=K\oplus L(K^{\perp})$. Let P denote the orthogonal projection from 
$C^{0,\alpha}(S_+^n)$
onto $K$, and $T:K \to \mathbb{R}^n$ be the isomorphism sending $x^i|_{S_+^n}$ to $e_i$, the $i$th coordinate basis. Define
$\tilde{P}=T\circ P$, that is, 
$$
 \tilde{P}(f)=\frac{2}{w_{n+1}}\left( \int_{S_+^n}f\, x^i\right) e_i 
$$
because
$$
\int_{S_+^n}x^i\, x^j=\frac{w_{n+1}}{2}\,\delta_{ij},
$$
where $w_{n+1}=\mathrm{Vol}(\mathbb{B}_1)$.

\begin{lemma}
\label{lemma4}
We have
\begin{equation}
\label{P(H)}
\tilde{P}\left(H(r,\tau ,0)\right)=-\frac{2\,w_n\,r^2}{(n+2)w_{n+1}}\,h^{\tau}_{jj,\underline{i}}\,e_i \, +\, O(r^3).
\end{equation}
\end{lemma}
{\it Proof}. From Lemma \ref{lemma2} and the fact that $\tilde{P}(x^{n+1})=\tilde{P}(x^ix^jx^{n+1})=0$  we have
$$
\tilde{P}\left(H(r,\tau ,0)\right)=\left[
-(n+4)h^{\tau}_{ij,\underline{k}}\,\tilde{P}\left(x^{n+1}x^ix^jx^k\right)
+2h^{\tau}_{ji,\underline{j}} \,\tilde{P}\left(x^{n+1}x^i\right)
\right]r^2+O(r^3)
$$
where
$$
\tilde{P}\left(x^{n+1}x^ix^jx^k\right)=\frac{2\left(\delta_{ij}\delta_{kl}+\delta_{ik}\delta_{jl}
+\delta_{il}\delta_{jk}\right)w_n}{(n+2)(n+4)w_{n+1}}\;e_l
$$
and
$$
\tilde{P}\left(x^{n+1}x^i\right) =\frac{w_n\,\delta_{ij}}{w_{n+1}(n+2)}\;e_j \, .
$$
Hence (\ref{P(H)}) follows. \hfill \fbox 

\begin{lemma}
\label{lemma5}
If $\varphi^{\tau}:\mathbb{S}_+^n\to \mathbb{R}$ is the solution of the Neumann problem
\begin{equation}
\label{Nprob}
\left\{
\begin{array}{l}
-(\Delta_{\mathbb{S}_+^n} \varphi +n\,\varphi)= h_{ii}^{\tau}\,x^{n+1}
-(n+3)\, h_{ij}^{\tau}\,x^ix^jx^{n+1} \\ \\
\displaystyle\frac{\partial \varphi}{\partial e_{n+1}}=0 \;\;on\; \partial(\mathbb{S}_+^n )
\end{array}
\right.
\end{equation}
then
\begin{equation}
\label{P(H)2}
\tilde{P}\left(H_{\varphi r}(0,\tau ,0)\varphi_{\tau}\right)=0
\end{equation}
and
\begin{equation}
\label{P(H)3}
\tilde{P}\left(H_{\varphi \varphi}(0,\tau ,0)\varphi_{\tau} \varphi_{\tau}\right)=0.
\end{equation}
for $\|\tau \|<r_p/2$.
\end{lemma}
{\it Proof}. The function $\varphi^{\tau}=\frac{1}{2}h_{ij}^{\tau}\,x^ix^jx^{n+1}$ is the solution of the
problem (\ref{Nprob}). Now we can use (\ref{H_rphi}), (\ref{H_phiphi}) and the fact 
$$
\tilde{P}(x^{i_{1}} \dots x^{i_{2k}}x^{n+1})=0,
$$
for all interge $k\geq 0$, to prove (\ref{P(H)2}) and (\ref{P(H)3}).

\section{Main Theorem}
\begin{defi} Consider $p\in \partial M$ and let $U$ be a neighborhood of $p$ on $M\cup \partial M$. A smooth codimension 1 foliation  $\mathscr{F}$ of $U\backslash \{p\}$ for a neighborhood $U$ of $p$ is called a free boundary
foliation centered at $p$, provided that its leaves are all closed and free boundary.
\end{defi}

\begin{teo}
\label{theorem1}
If $p\in \partial M$ is a nondegenerate critical point of the mean curvature function of $\partial M$, 
then there exist $\delta >0$ and smooth functions $\tau =\tau (r)$ and $\varphi =\varphi (r)$ with 
$\tau (0)=0$ such that $H(r,\tau (r), r \varphi (r))\equiv  n$ for all $0\leq r<\epsilon $. Hence the family 
$\mathscr{F}=\{S_r:=S_{r,\tau(r),r\varphi (r)}:\, 0\leq r<\epsilon \}$ is a smooth family of constant mean curvature spheres with $S_r$ having mean curvature $n/r$. Furthermore
$\mathscr{F}$ is a free boundary foliation centered at p.
\end{teo}
{\it Proof}. We will use the Taylor's formula with integral remainder 
$$
H(r,\tau , r\varphi )=n+\left[ H_{\varphi}(0,\tau ,0)\varphi +H_r(0,\tau ,0)\right] r
$$
$$
+\left[\frac{1}{2}\,H_{\varphi \varphi}(0,\tau ,0)\varphi \varphi +H_{\varphi r}(0,\tau ,0)\varphi +\frac{1}{2}\,H_{rr}(0,\tau ,0)\right]r^2+R(r,\tau ,\varphi)\,r^3
$$
where
$$
R(r,\tau ,\varphi)=\int_0^1(1-\eta )H_{\varphi r r}(\eta r,\tau ,0)\varphi \,d\eta
+\frac{1}{2}\int_0^1H_{\varphi \varphi r}(\eta r,\tau , 0 )\varphi \varphi \,d\eta
$$
$$
+\frac{1}{2}\int_0^1(1-\eta )^2H_{\varphi \varphi \varphi}(r,\tau , \eta \varphi )\varphi \varphi \varphi \,d\eta .
$$
We are interested in solving the equation $H(r,\tau , r\varphi )=n$, but first we are going to treat the equation
\begin{equation}
\label{PperpH=0}
P^{\perp}(H(r,\tau ,r\varphi)-n)=0,
\end{equation}
where $P^{\perp}$ denotes the $L^2$ orthogonal projection from $C_N^{0,\alpha}$ onto $K^{\perp}$. 

 By (\ref{H_phi}) and the fact $PL =PH_r(0,\tau ,0)=0$ we can write the equation in (\ref{PperpH=0}) 
 as follows (after division by $r$)
$$
L\varphi+H_r(0,\tau ,0)+\bar{R}(r,\tau ,\varphi)\,r=0,
$$
where
$$
\bar{R}(r,\tau ,\varphi)=
\frac{1}{2}\,P^{\perp}(H_{\varphi \varphi}(0,\tau ,0)\varphi \varphi )
+P^{\perp}(H_{\varphi r}(0,\tau ,0)\varphi )+P^{\perp}(\frac{1}{2}\,H_{rr}(0,\tau ,0))
$$
$$
+P^{\perp}(R(r,\tau ,\varphi))\,r.\hspace{5.5cm}
$$
Consider the mapping $G:[0,r_p/8)\times \mathbb{B}_{r_p/4} \times \mathbb{B}_{\delta_0} \to K^{\perp}$ given by
$$
G(r,\tau ,\varphi )=L\varphi+H_r(0,\tau ,0)+\bar{R}(r,\tau ,\varphi)\,r ,
$$ 
where $\mathbb{B}_{\delta_0}=\{\varphi \in K^{\perp}\,;\;\|\varphi \|_{C^{2,\alpha}}<\delta_0  \}$.

\vspace{0.2cm}
For $\tau =0$, let $\varphi_0 \in C_N^{2,\alpha}(\mathbb{S}^n)$ be a solution of the equation 
$$
L\varphi_0+H_r(0,0,0)=0.
$$
One sees that $G(0,0,\varphi_0 )=0$ and $G_{\varphi}(0,0,\varphi_0)=\left(-\Delta -n\right):K^{\perp}\to L(K^{\perp})$
is a bounded invertible linear transformation. By the implicit function theorem we can solve 
$P^{\perp}(H(r,\tau , r\,\varphi (r,\tau ))-n)=0$ for a function 
$ \varphi :[0,\delta )\times \mathbb{B}_{\delta}\to K^{\perp}$, 
for some $0<\delta \leq r_p/8$, with $ \varphi (0,0)=\varphi_0$. 
Furthermore
$$
\varphi_r(0,0)=-G_{\varphi}(0,0,\varphi_0)^{-1}G_r(0,0,\varphi_0) \hspace{6cm}
$$
$$
\hspace{1.4cm} =-\left(-\Delta -n\right)^{-1}\frac{\partial}{\partial r}(\bar{R}(r,\tau ,\varphi)\,r)|_{r=0,\varphi =
\varphi_0}
=-\left(-\Delta -n\right)^{-1}\bar{R}(0,0 ,\varphi_0)
$$
where 
$$
\bar{R}(0,0 ,\varphi_0)=
\frac{1}{2}\,H_{\varphi \varphi}(0,0,0)\varphi_0 \varphi_0 
+H_{\varphi r}(0,0,0)\varphi_0 .
$$
Since  
$$
L\varphi(r,\tau )+H_r(0,\tau ,0)+O(r)=0,
$$
we have, for $r=0$,
$$
L\varphi(0,\tau )+H_r(0,\tau ,0)=0.
$$
Then, by Lemma \ref{lemma5}, 
$$
\tilde{P}\left(H_{\varphi r}(0,\tau ,0)\varphi(0,\tau )\right)=
\tilde{P}\left(H_{\varphi \varphi}(0,\tau ,0)\varphi(0,\tau ) \varphi(0,\tau )\right)=0.
$$ 
On the other hand,
$$
\varphi(r,\tau )=\varphi(0,\tau )+r\int_0^1\varphi_r(\eta r,\tau )\,d\eta \, ,
$$
so that
\begin{equation}
\label{2eq}
\tilde{P}\left(H_{\varphi r}(0,\tau ,0)\varphi (r, \tau)\right)=r\, 
\tilde{P}\left(H_{\varphi r}(0,\tau ,0)\left(\int_0^1\varphi_r(\eta r,\tau )\,d\eta\right)\right)
\end{equation}
and
$$
\tilde{P}\left(H_{\varphi \varphi}(0,\tau ,0)\varphi (r, \tau) \varphi (r, \tau)\right)=
$$
$$
r\, 
\tilde{P}\left(H_{\varphi \varphi}(0,\tau ,0)\varphi (r, \tau)\left(\int_0^1\varphi_r(\eta r,\tau )\,d\eta\right)\right)
+O(r^2).
$$

Now we consider the mapping $(r,\tau )\mapsto H(r, \tau , r\varphi (r,\tau ))-n$ whose values lie in $K$ by the construction of $\varphi (r,\tau )$.
Let us solve the equation 
$$
H(r, \tau , r\varphi (r,\tau ))-n=0,
$$ 
which is equivalent to equation $\tilde{P}(H(r,\tau , r\varphi (r,\tau ))-n)=0$ and, after division by $r^2$, it is equivalent to
$$
\frac{1}{2}\,\tilde{P}\left(H_{\varphi \varphi}(0,\tau ,0)\varphi \varphi\right) +
\tilde{P}\left(H_{\varphi r}(0,\tau ,0)\varphi\right) +
\frac{1}{2}\,\tilde{P}\left(H_{rr}(0,\tau ,0)\right)+\tilde{P}\left(R(r,\tau ,\varphi)\right)r=0,
$$
where $\varphi =\varphi (r,\tau )$. By (\ref{2eq}) and Lemma 0.4 the above equation may
be written as follows
$$
-\frac{2\,w_n}{(n+2)w_{n+1}}\,h^{\tau}_{jj,\underline{i}}\,e_i +R_1(r,\tau ,\varphi)r=0,
$$
with
$$
R_1(r,\tau ,\varphi)=\tilde{P}\left(R(r,\tau ,\varphi)\right)+r\, 
\tilde{P}\left(H_{\varphi r}(0,\tau ,0)\left(\int_0^1\varphi_r(\eta r,\tau )\,d\eta\right)\right)
$$
$$
\hspace{1.5cm}+\,r\, 
\tilde{P}\left(H_{\varphi \varphi}(0,\tau ,0)\varphi (r, \tau)\left(\int_0^1\varphi_r(\eta r,\tau )\,d\eta\right)\right)
+O(r^2).
$$
In order to solve the above equation, consider $F:[0,\delta)\times B(0,\delta)\to \mathbb{R}^n$ defined by
\begin{equation}
\label{Frtau}
F(r,\tau )=-\frac{2\,w_n}{(n+2)w_{n+1}}\,h^{\tau}_{jj,\underline{i}}\,e_i +R_1(r,\tau ,\varphi)\, r.
\end{equation}
By the assumption $h^{\tau}_{,\underline{i}}|_{\tau =0}=0$, we have $F(0,0)=0$
and the Hessian matrix
$$
\frac{\partial F}{\partial \tau}(0,0)=\left(\frac{\partial}{\partial \tau_j}h^{\tau}_{,\underline{i}}\bigg\vert_{\tau=0}\right)_{i,j}
$$
is nonsingular. Applying the implicit function theorem we obtain a solution $\tau =\tau(r)$
of the equation $H(r,\tau ,r\varphi(r,\tau))=n$ around $(r,\tau)=(0,0)$, $0\leq r <\epsilon$,
for some $0<\epsilon \leq \delta$. 

\subsection{The Foliation}
It is clear that $\mathscr{F}=\{S_r=S_{r,\tau(r),r\varphi(r)}\, ; \;\;0<r\leq \epsilon \}$
is a smooth family of embedded constant mean curvature hemisphere with $S_r$ having mean curvature $n/r$. We need to show that this family constitutes a foliation.

In order to prove this we consider the aplication $X:(0,\epsilon ) \times \mathbb{B}_1 \to M$ given by
$$
X(r,x)=\varphi^{\tau(r)}(r(1+r\varphi(r)(x,t(x)))(x,t(x)))
$$
where $t(x)=\sqrt{1-|x|^2}$ and $\varphi^{\tau}$ is the Fermi coordinate system defined in (\ref{varphi^tau}).
Observe that $X(B_1 ,r)=S_r$. 
Thus, it is sufficient to prove that $X$ is a parametrization of $M$, for small $\epsilon >0$.
It is enough to prove that
$$\Upsilon(r,x)=(\varphi^0)^{-1}(X(r,x))$$
is an immersion, where $\varphi^0$ is a Fermi coordinate system centred at $p$ defined in (\ref{varphi0}).

\noindent{\it Claim:} The function $\tau =\tau (r)$ satisfies $\tau (r)=O(r^2)$. 

We have that $\tau (r)$ is a solution of the equation $F(r,\tau (r))=0$, where $F$ is defined in (\ref{Frtau}).
By the implicit function theorem
$$
\frac{\partial \tau}{\partial r}(0)=-\left(\frac{\partial f_i}{\partial \tau_j}(0)\right)^{-1}\left(\frac{\partial F_1}{\partial r}(0),
\dots ,\frac{\partial F_n}{\partial r}(0)\right)
$$
and
$$
\frac{\partial F_i}{\partial r}(0,0)=\langle R_1(0,0 ,\varphi_0), e_i\rangle
=\langle \tilde{P}(R(0,0 ,\varphi_0)), e_i\rangle =0.
$$
Then,
\begin{equation}
\label{tau_r}
\frac{\partial \tau}{\partial r}(0)=0.
\end{equation}

Now
$$
\frac{\partial \Upsilon}{\partial r}(r,x)=\left(d_{x}\Upsilon \right)\left((1+r\varphi(r)(x,t))(x,t)+r((1+r\varphi(r)(x,t))_r(x,t)\right)
$$
$$
+\left(\frac{\partial \Upsilon}{\partial r}\right)(r(1+r\varphi(r)(x,t))(x,t))
$$
and
$$
\frac{\partial \Upsilon}{\partial r}(0,0)=\frac{\partial}{\partial \tau^i}\left( (\varphi^0)^{-1}\circ \varphi^{\tau (r)} \right)\bigg\vert_{\tau =0}
\frac{\partial \tau^i}{\partial r}(0)=0 .
$$
Using (\ref{tau_r}) we conclude that 
$$
(\partial \Upsilon /\partial r)(0,x)=(x, t(x)),
$$
 i.e., 
 $$
 \Upsilon (r,x)=r(x, t(x))+O(r^2).
$$ 
Consequently 
$$
(\partial \Upsilon /\partial r)(r,x)=(x, t(x))+O(r),
$$ 
$$
(\partial \Upsilon /\partial x_i)(r,x)=re_i-r(x_i/t(x))e_t+O(r^2)
$$ 
and
$$
\mathrm{det}\left(\frac{\partial \Upsilon}{\partial r}\;\; \frac{\partial \Upsilon}{\partial x_1}\;\dots \;\frac{\partial \Upsilon}{\partial x_n}\right)(r,x)=r^2\left(\frac{1}{t(x)}+O(r)\right)>0,
$$
for all $0<r<\epsilon '\leq \epsilon$ and some $\epsilon '$ enought small.

\end{document}